\documentclass[12pt,a4]{article}

\usepackage[left=2.5cm,right=2.5cm,top=2.5cm,bottom=2.5cm,a4paper]{geometry}

\usepackage[leqno]{amsmath}
\usepackage[all]{xy}
\usepackage[mathscr]{euscript}
\usepackage{amssymb}
\usepackage{amscd}

\newtheorem{lemma}{Lemma}[section]
\newtheorem{theorem}[lemma]{Theorem}
\newtheorem{prop}[lemma]{Proposition}

\newtheorem{remark}[lemma]{Remark}

\newcommand{\pf}{\noindent{\em Proof: }}
\newcommand{\epf}{\hfill\hbox{\rule{3pt}{6pt}}\\}

\newcommand{\forme}[1]{}
\begin{document}

\title{The distance-regular graphs such that  all of its second largest local eigenvalues are at most one}

\author{Jack H. Koolen,
 Hyonju Yu\\
{\small {\tt koolen@postech.ac.kr} ~~
{\tt lojs4110@postech.ac.kr}}\\
{\footnotesize{Department of Mathematics,  POSTECH, Pohang 790-785, South Korea}}}

\date{\today}

\maketitle

{\bf Dedicated to Professor Dragos Cvetkovi\'c on the occasion of his 70th birthday}\\
\\

\begin{abstract}
In this paper, we classify distance regular graphs such that all of its second largest local eigenvalues are at most one. Also we discuss the consequences for the smallest eigenvalue of a distance-regular graph.
These extend  a result by the first author, who classified the distance-regular graph with smallest eigenvalue $-1-\frac{b_1}{2}$. \end{abstract}

\section{Introduction}
Koolen \cite{Koolen} classified the  distance-regular graphs with smallest eigenvalue $-1-\frac{b_1}{2}$.

\begin{theorem}(\cite{Koolen})\label{koo}
Let $\Gamma$ be a distance-regular graph with diameter $D$ and smallest eigenvalue $-1-\frac{b_1}{2}$. Then either $a_1\leq 1$ or one of the following holds:
\begin{itemize}
\item[(I)] $D=2$ and
\begin{itemize}
\item[(a)] $\Gamma$ is a complete multipartite graph $K_{n\times t}$ with $n\geq 4$, $t\geq 2$;
\item[(b)] $\Gamma$ is the complement of  an $n\times n$ grid with $n\geq 4$;
\item[(c)] $\Gamma$ is the complement of a triangular graph $T(n)$, with $n\geq 5$;
\item[(d)] $\Gamma$ is the complement of the Petersen graph;
\item[(e)] $\Gamma$ is the complement of the Shrikhande graph;
\item[(f)] $\Gamma$ is the complement of one of the three Chang graphs;
\end{itemize}
\item[(II)] $D=3$ and
\begin{itemize}
%\item[(a)] $\Gamma$ is the icosahedron with intersection array $\{5 ,2 ,1 ;1 ,2 ,5 \}$;
\item[(a)] $\Gamma$ is the Johnson graph $J(6,3)$ with intersection array $\{9, 4, 1; 1, 4, 9 \}$;
\item[(b)] $\Gamma$ is the distance-2 graph of the halved 6-cube with intersection array $\{15, 8, 1; 1, 8, 15 \}$;
\item[(c)] $\Gamma$ is the Gosset graph with intersection array $\{27, 16, 1; 1, 16, 27 \}$;
\end{itemize}
\item[(III)] $D=4$ and $\Gamma$ is the Conway-Smith graph with intersection array $\{10, 6, 4, 1; 1, 2, 6, 10\}$.
\end{itemize}
\end{theorem}

Note that Case II of Theorem 1.1 was forgotten in \cite{Koolen}.
\\

For a vertex $x$ of a graph $\Gamma$, let $\Delta(x)$ be the local graph of $\Gamma$, i.e. the subgraph induced on the neighbours of $x$.  In this paper, we extend this Theorem \ref{koo} as follows. We determine the distance-regular  graphs such that for all $x$ the second largest eigenvalue of $\Delta(x)$ is at most one. Our main result is:
\begin{theorem}\label{main}
Let $\Gamma$ be a distance-regular graph with diameter $D$ such that for all vertices $x$, the second largest eigenvalue of the local graph of $x$, $\Delta(x)$, is at most one. Then, either $a_1\leq 1$ or, one of the following holds:
\begin{itemize}
\item[(I)] $D=1$ and $\Gamma$ is the complete graph $K_n$ with $n\geq 4$;
\item[(II)] $D=2$ and
\begin{itemize}
\item[(a)] $\Gamma$ is a complete multipartite graph $K_{n\times t}$ with $n\geq 4$, $t\geq 2$;
\item[(b)] $\Gamma$ is a complement of an $n\times n$ grid with $n\geq 4$;
\item[(c)] $\Gamma$ is the complement of a  triangular graph $T(n)$, with $n\geq 5$;
\item[(d)] $\Gamma$ is the complement of the Petersen graph;
\item[(e)] $\Gamma$ is the complement of the Shrikhande graph;
\item[(f)] $\Gamma$ is the complement of one of the three Chang graphs;
\item[(g)] $\Gamma$ is the Shrikhande graph;
\item[(h)] $\Gamma$ is the Clebsch graph;
\item[(i)] $\Gamma$ is the Paley graph with 13 vertices;
\item[(j)] $\Gamma$ is the Paley graph with 17 vertices;
\item[(k)] $\Gamma$ has possibly $\{12, 6; 1, 6\}$,
%$\{15, 6; 1, 9\}$,
$\{15, 8; 1, 6\}$, $\{18, 10; 1, 6\}$, $\{21 ,12; 1, 6\}$, $\{21, 12; 1, 9\}$ or $\{27, 16; 1, 12\}$ as its intersection array;
\end{itemize}
\item[(III)] $D=3$ and
\begin{itemize}
\item[(a)] $\Gamma$ is the icosahedron with intersection array $\{5, 2, 1; 1, 2, 5 \}$;
\item[(b)] $\Gamma$ is the Johnson graph $J(6,3)$ with intersection array $\{9, 4, 1; 1, 4, 9 \}$;
\item[(c)] $\Gamma$ is the Doro graph with intersection array $\{10, 6, 4; 1, 2, 5 \}$;
\item[(d)] $\Gamma$ is the distance 2-graph of the halved 6-cube with intersection array $\{15, 8, 1; 1, 8, 15 \}$;
\item[(e)] $\Gamma$ is the unique locally folded 5-cube distance-regular graph with intersection array\\ $\{16 ,10 ,1 ;1 ,5 ,16 \}$;
\item[(f)] $\Gamma$ is the Gosset graph with intersection array $\{27, 16, 1; 1, 16, 27 \}$;
\end{itemize}
\item[(IV)] $D=4$ and $\Gamma$ is the Conway-Smith graph with intersection array $\{10, 6, 4, 1; 1, 2, 6, 10\}$.
\end{itemize}

\end{theorem}

As a consequence of Theorem \ref{main} and some results of Terwilliger and Hoffman, we can also extend Theorem \ref{koo} (and \cite[Theorem 4.4.3]{bcn}) in the following way.

\begin{theorem}\label{main2}
Let $0<\alpha < 1 + \sqrt{2}$. Then there exists $K = K(\alpha)$, such that any distance-regular graph with diameter $D$ at least three, valency $k \geq K$, $a_1 \geq 2$ and smallest eigenvalue at most $-1 -\frac{b_1}{\alpha}$ is one of the following graphs:
\begin{itemize}
\item[(I)] $D=3$ and
\begin{itemize}
\item[(a)] $\Gamma$ is the icosahedron with intersection array $\{5, 2, 1; 1, 2, 5 \}$;
\item[(b)] $\Gamma$ is the Johnson graph $J(6,3)$ with intersection array $\{9, 4, 1; 1, 4, 9 \}$;
%\item[(c)] $\Gamma$ is the Doro graph with intersection array $\{10, 6, 4; 1, 2, 5 \}$;
\item[(c)] $\Gamma$ is the distance-2 graph of the  halved 6-cube with intersection array $\{15, 8, 1; 1, 8, 15 \}$;
%\item[(e)] $\Gamma$ is the unique locally folded 5-cube distance-regular graph with intersection array $\{16, 10, 1; 1, 5, 16 \}$;
\item[(d)] $\Gamma$ is the Gosset graph with intersection array $\{27, 16, 1; 1, 16, 27 \}$;
\end{itemize}
\item[(IV)] $D=4$ and $\Gamma$ is the Conway-Smith graph with intersection array $\{10, 6, 4, 1; 1, 2, 6, 10\}$.
\end{itemize}
\end{theorem}

This paper is organized as follows: in the next section we will give definitions and preliminaries, in Section 3, we give some results that we will use in the proofs of Theorems \ref{main} and \ref{main2}.
In Section 4 we will give the proof of Theorem \ref{main}, in the first part we consider diameter at least three and the second part we consider diameter two. In the last section, we show Theorem \ref{main2}.

\section{Definitions and preliminaries}
All the graphs considered in this paper are finite, undirected and
simple (for unexplained terminology and more details, see \cite{bcn}). Suppose
that $\Gamma$ is a connected graph with vertex set $V(\Gamma)$ and edge set $E(\Gamma)$, where $E(\Gamma)$ consists of unordered pairs of two adjacent vertices. The distance $d(x,y)$ between
any two vertices $x,y$ of $\Gamma$
is the length of a shortest path connecting $x$ and $y$ in $\Gamma$.
We denote by $\overline{\Gamma}$ the complement of $\Gamma$.

Let $\Gamma$ be a connected graph. For a vertex $x \in V(\Gamma)$, define $\Gamma_i(x)$ as the set of
vertices which are at distance precisely $i$ from $x~(0\le i\le
D)$, where $D:=\max\{d(x,y)\mid x,y\in V(\Gamma)\}$ is the diameter
of $\Gamma$. In addition, define $\Gamma_{-1}(x) = \Gamma_{D+1}(x)
= \emptyset$. We write $\Gamma(x)$ instead of $ \Gamma_1(x)$. The adjacency matrix $A$ of graph $\Gamma$ is the (0,1)-matrix whose rows and columns are indexed by the vertex set $V(\Gamma)$ and the $(x,y)$-entry is $1$ whenever $x$ and $y$ are adjacent (denoted by $x \sim y$) and 0 otherwise.
The eigenvalues (respectively, the spectrum) of the graph $\Gamma$ are the eigenvalues (respectively, the spectrum) of  $A$. We denote the second largest eigenvalue of respectively a graph $\Gamma$ or a square matrix $Q$ with only real eigenvalues by $\theta_1(\Gamma)$, respectively $\theta_1(Q)$.

For a connected graph $\Gamma$, the {\em local graph} $\Delta(x)$ of a vertex $x\in V(\Gamma)$ is the subgraph induced on $\Gamma(x)$ in $\Gamma$.
%  and the induced subgraph on $\Gamma(x)\cup\{x\}$ is denoted by $\widehat{\Delta(x)}$.

For a graph $\Gamma$, a partition $\Pi= \{P_1, P_2, \ldots ,
P_{\ell}\}$ of the vertex set $V(\Gamma)$ is called {\em
equitable} if there are constants $\beta_{ij}$ such that each vertex $x \in P_i$ has exactly $\beta_{ij}$ neighbors in $P_j$ ($1\leq i, j \leq \ell$). The {\em quotient
matrix} $Q(\Pi)$ associated with the equitable partition $\Pi$ is
the $\ell \times \ell$ matrix whose $(i,j)$-entry equals $\beta_{ij}$ ($1\leq i,j\leq \ell$). Note that the eigenvalues of the quotient matrix $Q(\Pi)$ are also eigenvalues (of the adjacency matrix $A$) of $\Gamma$.

A connected graph $\Gamma$ with diameter $D$ is called {\em{distance-regular}} if there are integers $b_i, c_{i}$ $(0\leq i\leq D)$ (where $b_D= 0 = c_0$) such that for any two vertices $x, y \in V(\Gamma)$ with $d(x, y)=i$, there are precisely $c_i$ neighbors of $y$ in $\Gamma_{i-1}(x)$ and $b_i$ neighbors of $y$ in $\Gamma_{i+1}(x)$. In particular, any distance-regular graph  is regular with valency $k := b_0$. Note that a {\em non-complete, connected  strongly regular graph} is just a distance-regular graph with diameter two. In this case we say that $\theta_1$ and $\theta_2$ are the non-trivial eigenvalues.
We define $a_i := k-b_i-c_i \ \ (1 \leq i \leq D)$ for notational convenience.  Note that $a_i=\mid\Gamma(y)\cap\Gamma_i(x)\mid$ holds for any two vertices $x, y$ with $d(x, y)=i$ $(1\leq i\leq D).$ For a distance-regular graph $\Gamma$ and a vertex $x\in V(\Gamma)$, we denote $k_i:=|\Gamma_i(x)|$. It is easy to see that $k_i = \frac{b_0 b_1 \cdots b_{i-1}}{c_1 c_2 \cdots c_i}$ and hence  $k_i$ does not depend on the vertex $x$.  The numbers $a_i$, $b_{i-1}$ and $c_i$ $(1\leq i\leq D)$ are called the {\em{intersection~numbers}} of the distance-regular graph $\Gamma$, and the array $\{b_0,b_1,\ldots,b_{D-1};c_1,c_2,\ldots,c_D\}$ is called the {\em{intersection~array}} of $\Gamma$.

The next lemma gives some elementary properties concerning the intersection numbers.

\begin{lemma}(Cf. \cite[Proposition 4.1.6]{bcn})\label{pre}{\ \\}
Let $\Gamma$ be a distance-regular graph with valency $k$ and diameter $D$. Then the following holds:\\
(1) $k=b_0> b_1\geq \cdots \geq b_{D-1}~;$\\
(2) $1=c_1\leq c_2\leq \cdots \leq c_{D}~;$\\
(3) $b_i\ge c_j$ \mbox{ if } $i+j\le D~.$
\end{lemma}

We will refer to the following theorem as the interlacing theorem.

\begin{theorem}\label{0}(Cf.  \cite[Theorem 9.1.1]{god})  Let $m \geq n$ be two positive integers. Let $A$ be an  $n\times n$ matrix, which is similar to a (real) symmetric matrix, and let $B$ be a principal $m \times m$ submatrix of $A$. Then, for $i=1,\ldots , m$, $$\theta_{n-m+i}(A)\leq \theta_i(B)\leq \theta_i(A)$$
holds, where $A$ has eigenvalues $\theta_1(A) \geq \theta_2(A) \geq \cdots \geq \theta_n(A)$ and $B$ has eigenvalues  $\theta_1(B) \geq \theta_2(B) \geq \cdots \geq \theta_m(B)$.
\end{theorem}

For $\Pi= \{P_1, P_2, \ldots, P_t\}$ a partition of the vertex set of a graph $\Gamma$ the {\em quotient matrix} with respect to $\Pi$ is the $t \times t$-matrix $B$ whose $(i,j)$-entry equals $\sum_{x \in P_i} $ $\frac{(\# \mbox{ of neighbours of $x$ in $P_j$)}}{|P_i|}$.

We also need the following interlacing result.

\begin{theorem}\label{god}(Cf. \cite[Lemma 9.6.1]{god}) Let $\Gamma$ be a graph with $\nu$ vertices and eigenvalues $\theta_0 \geq \theta_1 \geq \cdots \geq \theta_{\nu-1}$. Let $\Pi= \{P_1, P_2, \ldots, P_t\}$ be a partition of the vertex set  of $\Gamma$ with quotient matrix $B$ with respect to $\Pi$  Then, for $i=1,\ldots , m$, $$\theta_{\nu-t+i}\leq \theta_i(B)\leq \theta_i$$
holds, where $B$ has eigenvalues  $\theta_1(B) \geq \theta_2(B) \geq \cdots \geq \theta_t(B)$. Moreover if the interlacing is tight (there exists $\ell$ such that $\theta_i(B) = \theta_i$ for $i \leq \ell$ and
$\theta_i(B) = \theta_{\nu - t +i}$ for $i > \ell$), then $\Pi$ is equitable.
\end{theorem}

The following theorem summarizes some elementary results on strongly regular graphs.

\begin{theorem}\label{srg}(Cf. \cite[Theorem 1.3.1 and Proposition 1.3.2]{bcn}
Let $\Gamma$ be a non-complete connected strongly regular graph with $\nu$ vertices, distinct eigenvalues $k>\theta_1>\theta_2$ and intersection numbers $c_2$ and $b_1$. Then the followings hold;
\begin{itemize}
\item[(i)] $k=c_2-\theta_1\theta_2$;
\item[(ii)] the multiplicity of $\theta_1$ equals $\frac{(\theta_2+1)k(k-\theta_2)}{c_2(\theta_2-\theta_1)}$;
\item[(iii)] $(\theta_1 +1)(\theta_2 +1) = -b_1$;
\item[(iv)] If $\Gamma$ is not a conference graph, then $\theta_1$ and $\theta_2$ are integers;
%\begin{claim}\label{claim}
\item[(v)] If $\Gamma$ contains a coclique $C$ of size $\gamma$, then $\theta_1\leq \frac{(\nu-\gamma)(k-c_2)}{\gamma k}$.
%\pf
%Since $\Gamma$ has a coclique $C$ of size $t$, we can consider the quotient matrix $\left[
%\begin{array}{cc}
% 0  & k \\
% \alpha & k-\alpha
%\end{array} \right]$ of $\{VC, V\Gamma \setminus VC\}$ where $\alpha=\frac{tk}{\nu-t}$. Then $-\theta_2 \geq \alpha$. So, we obtain that $a_2=\theta_1\theta_2\geq \frac{tk}{\nu-t}\theta_1$.
%\epf
\end{itemize}
 \end{theorem}

The next theorem summarizes the results on regular graphs with smallest eigenvalue at least $-2$ and is in essence due to Cameron et al. \cite{CGSS}.

\begin{theorem}(Cf. \cite[Proposition 3.12.2]{bcn})\label{cla1}
Let $\Gamma$ be a connected regular graph with $\nu$ vertices, valency $k$, and smallest eigenvalue at least $-2$. Then one of the following holds:
\begin{itemize}
\item[(i)] $\Gamma$ is the line graph of a regular connected graph;
\item[(ii)] $\Gamma$ is the line graph of bipartite semiregular connected graph;
\item[(iii)] $\nu=2(k+2)\leq 28$ and $\Gamma$ is an induced subgraph of $E_7(1)$;
\item[(iv)] $\nu=\frac{3}{2}(k+2)\leq 27$ and $\Gamma$ is an induced subgraph of Schl\"{a}fli graph;
\item[(v)] $\nu=\frac{4}{3}(k+2)\leq 16$ and $\Gamma$ is an induced subgraph of Clebsch graph;
\item[(vi)] $\nu=k+2$ and $\Gamma$ is a $ K_{m\times 2}$ for some $m\geq3$.
\end{itemize}
\end{theorem}

\begin{remark}
(i) There are 187 regular connected graphs with smallest eigenvalue at least $-2$, which are not line graphs, see for example \cite[p. 91]{BCS}.
\\
(ii) In Cases (iii)-(v) of Theorem \ref{cla1}, one can say more by inspecting the 187 regular graphs of (i), namely there are only 5 graphs (all of which have 22 vertices) which are not an induced subgraph of one of the Schl\"{a}fli graph or
the three Chang graphs.
\end{remark}

The following result was originally shown by J.J. Seidel \cite{Sei}.

\begin{theorem}(Cf. \cite[Proposition 3.12.4]{bcn})\label{cla2}
Let $\Gamma$ be a connected strongly regular graph with smallest eigenvalue $-2$. Then $\Gamma$ is a triangular graph $T(n) \ (n\geq 5)$, a square grid $n\times n \ (n\geq 3)$, a complete multipartite graph $K_{n\times 2} \ (n\geq 2)$, or one of the graphs of Petersen, Clebsch, Schl\"{a}fli, Shrikhande, or Chang.
\end{theorem}

Terwilliger \cite{Ter} showed the following diameter bound for distance-regular graphs containing an induced quadrangle.
\begin{theorem}\label{Ter}(Cf. \cite[Theorem 5.2.1 and Corollary 5.2.2]{bcn})
Let $\Gamma$ be a distance-regular graph with diameter $D$. If $\Gamma$ contains an induced quadrangle, then $$ c_i-b_i\geq c_{i-1}-b_{i-1}+a_1+2 \ \ \ \ (i=2, \cdots, D).$$
In particular, $$D \leq \frac{k+c_D}{a_1 +2} \leq \frac{2k}{a_1 +2}.$$
\end{theorem}

Terwilliger \cite{Ter1} also determined the graphs which reach this diameter bound.

\begin{theorem}\label{Tercla}(Cf. \cite[Theorem 5.2.3]{bcn})
Let $\Gamma$ be a distance-regular graph with diameter $D$ at least $\frac{k+c_D}{a_1+2}$. Then one of the following holds:
\begin{itemize}
\item[(i)] $\Gamma$ is a Terwilliger graph;
\item[(ii)] $\Gamma$ is a strongly regular graph with smallest eigenvalue $-2$;
\item[(iii)] $\Gamma$ is a Hamming graph, a Doob graph, a Johnson graph, a halved cube, or the Gosset graph.
\end{itemize}

\end{theorem}

Tewilliger \cite{Ter2} also showed the following result on the eigenvalues of a distance-regular graph.

\begin{theorem}\label{4.4.3}(Cf. \cite[Theorem 4.4.3]{bcn})
Let $\Gamma$ be a distance-regular graph with diameter $D$ at least three and  distinct eigenvalues $k=\theta_0>\theta_1> \cdots >\theta_D$.
Let $x$ be a vertex of $\Gamma$ and let $\Delta(x)$ have eigenvalues $a_1 =\lambda_1 \geq \lambda_2 \geq \cdots \geq \lambda_k$.
Then $-1 - \frac{b_1}{\theta_D +1} \geq \lambda_2 \geq \lambda_k \geq -1 - \frac{b_1}{\theta_1 +1}$.
\end{theorem}

%J.I. Hall \cite{hall} has classified the locally Petersen graphs.

%\begin{theorem}(Cf. \cite[Theorem 1.16.5]{bcn})
%There are up to isomorphism three connected locally Petersen graphs, namely:
%\begin{itemize}
%\item[(i)] the complement of the triangular graph $T(7)$;
%\item[(ii)] the Conway-Smith graph with intersection array $\{10, 6, 4, 1; 1, 2, 6, 10\}$;
%\item [(iii)] the Doro graph with intersection array $\{10 ,6 ,4 ;1 ,2 ,5 \}$.
%\end{itemize}
%\end{theorem}

Recall that a Terwilliger distance-regular graph is a distance-regular graph such that the induced subgraph on the common neighbours of any two vertices at distance two is complete.

Koolen \cite{Koolen} showed:
\begin{prop}\label{ter}\cite[Proposition 9]{Koolen}
Let $\Gamma$ be a distance-regular Terwilliger graph with $c_2\geq 2$. If $\Gamma$ has a vertex $x$ such that $\theta_1(\Delta(x))\leq1$, then $\Gamma$ is one of the following:
\begin{itemize}
\item[(i)] the icosahedron with intersection array $\{5 ,2 ,1 ;1 ,2 ,5 \}$;
\item[(ii)] the Doro graph with intersection array $\{10 ,6 ,4 ;1 ,2 ,5 \}$;
\item[(iii)] the Conway-Smith graph with intersection array $\{10, 6, 4, 1; 1, 2, 6, 10\}$.
\end{itemize}
\end{prop}

The following result shows a construction of  antipodal distance-regular graphs with diameter three.

\begin{prop}(\cite[Proposition 12.5.3]{bcn})
Let $q = rm +1$ be a prime power, where $r >1$ is an integer and either $m$ is even or $q$ is a power of two. Let $V$ be a vector space of dimension two over the finite field with $q$ elements, GF$(q)$. Let $V$ be provided with a non-degenerate symplectic form $B$. Let $K$ be the subgroup of the multiplicative GF$(q)^* =$ GF$(q) \setminus \{0\}$ of index $r$, and let $b \in$ GF$(q)^*$. Then the graph $\Gamma$ with vertex set $\{ Kv \mid v \in V \setminus \{0\} \}$, where $\{ Ku, Kv\}$ is an edge
if $B(u,v) \in bK$ and $Ku \neq Kv$ is distance-regular with diameter three,  with $r(q+1)$ vertices and
intersection array $\{q, q-m-1, 1; 1, m, q\}$, and it is anantipodal $r$-cover of the complete graph $K_{q+1}$.
\end{prop}

\begin{remark} \label{example}
For $q=16$ and $r=3$ we obtain a locally folded 5-cube distance-regular graph with intersection array $\{16, 10, 1; 1, 5, 16\}$. This is the only known example of a distance-regular graph with this intersection array.
We will show in Proposition \ref{unique}  that there is a unique distance-regular graph with intersection array $\{ 16, 10, 1; 1, 5, 16\}$, that is locally the folded 5-cube, a result that also was obtained by \cite{boinck}, cf. \cite[p. 386]{bcn}.
\end{remark}

\section{Some useful results}

In this section we give some results which will be helpful to show our main results. First we give some sufficient conditions for a local graph of a distance-regular graph to be connected and coconnected, that is its complement is connected.

\begin{prop}\label{conn}
Let $t$ be a positive integer. Let $\Gamma$ be a distance-regular graph such that the local graph $\Delta(x)$ has second largest eigenvalue at most $t$ for any vertex $x$ of $\Gamma$. Then the following statements  hold:
\begin{itemize}
\item[(i)] If $a_1>t$, then for any vertex $x$, the local graph $\Delta(x)$ is connected;
\item[(ii)] If $\Gamma$ is not complete multipartite, then for any vertex $x$, the complement of $\Delta(x)$, $\overline{\Delta(x)}$, is connected.
\end{itemize}
\end{prop}
\pf
$(i)$: This follows immediately from Theorem \ref{0}. $(ii)$: If the complement of $\Delta(x)$ is not connected, then $\theta_{\min}(\triangle(x))=-1- b_1$. This means that $\theta_{\min}(\Gamma)\leq -1-b_1$, and hence, by \cite[Theorem 4.4.4]{bcn}, for all vertices $x$, the second largest eigenvalue of $\Delta(x)$ is at most $-1-\frac{b_1}{1+(-1-b_1)}=0$. So $\overline{\Delta(x)}$ is a disjoint union of cliques. This means that $\overline{\Delta(x)}$ is complete multipartite for all $x$ and hence $\Gamma$ is complete multipartite.
\epf

In the next lemma we show a lower bound for the intersection number $c_2$ for a distance-graph $\Gamma$ such that for some vertex $x$, the local graph $\Delta(x)$ is the complement of a line graph.

\begin{lemma}\label{c2}
Let $\Gamma$ be a distance-regular graph with valency $k$ such that for a vertex $x$ the complement of $\Delta(x)$ is the line graph of a graph $\Sigma$.
\\
(i) If $\Sigma$ is $t$-regular, with $t \geq 2$, then $c_2 \geq k - 3t+3$.\\
(ii) If $\Sigma$ is semiregular, with degrees $s, t$ satisfying $2 \leq s < t$, then $c_2 \geq k -  2s-t +3$.
\end{lemma}
\pf
(i) Take two distinct edges $uv$ and $uw$ of $\Sigma$. The number of edges that contain one of $\{u,v,w\}$ is at most $3t-4$. This shows (i).
\\
(ii) Similar argument.
\epf

As a consequence of Theorem \ref{cla2}, we have:

\begin{prop}\label{cla3}
\begin{itemize}
 \item [(i)] If $\Gamma$ is a strongly regular graph with $a_1\geq 2$ with smallest eigenvalue at least $-2$ and if for any vertex $x$ the eigenvalue $\theta_1(\Delta(x))$ is at most one, then $\Gamma$ is the Shrikhande graph or the Clebsch graph.
 \item [(ii)] If $\Gamma$ is a coconnected strongly regular graph with $a_1\geq 2$ such that $\theta_1\leq 1$, then $\Gamma$ is the complement of an $n\times n$ grid $(n\geq 4)$, a triangular graph $T(n) \ \ (n\geq 5)$, the Petersen graph, Shrikhande graph or a Chang graph.
 %\item [(i)] If $\Gamma$ is not coconnected strongly regular graph with $a_1\geq 2$ such that $\theta_1\leq 1$, then $\Gamma$ is complete multipartite graph. \\
 \end{itemize}
\end{prop}

%\begin{lemma}\label{disc}
%Let $\Gamma$ be a strongly regular graph with vertex $x$ and distinct eigenvalues $k>\theta_1>\theta_2$ such that the %complement of its local graph of $x$ is not connected. Then
%$$ k + \theta_1(k-1)\leq (-\theta_2)(1+k_2)$$
%\end{lemma}
%\pf
%Let $\theta_k'$ be the smallest eigenvalue of the local graph of $x$. Then $b_1=-\theta_k'-1$. By Theorem \ref{srg}, %$\theta_2\neq -1-b_1$ as $\theta_1\neq 0$. So multiplicity of $\theta_2$ (say $m_2$) is at most $1+k_2$ and this %implies that multiplicity of $\theta_1$ (say $m_1$) is at least $k-1$. Therefore $k+\theta_1(k-1)\leq k+\theta_1 %m_1=(-\theta_2)m_2\leq(-\theta_2)(1+k_2)$.
%\epf

In the next result, we classify the Taylor graphs such that for some vertex $x$ the local graph of $x$, $\Delta(x)$, has second largest eigenvalue at most one.

\begin{prop}\label{tay}
Let $\Gamma$ be a Taylor graph with $a_1\geq2$ such that for some vertex $x$, $\theta_1(\Delta(x))\leq1$. \\
Then $\Gamma$ is one of the following graphs:
\begin{itemize}
\item[(i)] the icosahedron with intersection array $\{5 ,2 ,1 ;1 ,2 ,5 \}$;
\item[(ii)] the Johnson graph $J(6,3)$ with intersection array $\{9 ,4 ,1 ;1 ,4 ,9 \}$;
\item[(iii)] the distance-2 graph of the halved 6-cube with intersection array $\{15 ,8 ,1 ;1 ,8 ,15 \}$;
\item[(iv)] the Gosset graph with intersection array $\{27 ,16 ,1 ;1 ,16 ,27 \}$.
\end{itemize}
\end{prop}

\pf
Let $\Gamma$ be a Taylor graph with $a_1\geq 2$ such that for some vertex $x$, $\theta_1(\Delta(x))\leq1$. Then by \cite[Theorem 1.5.3]{bcn}, it is known that $\Delta(x)$ is a strongly regular graph, say with parameters $(\overline{v}=k, \overline{k}=a_1, \overline{\lambda}, \overline{\mu})$, where $2\overline{\mu}=\overline{k}=a_1$ holds. As $\theta_1(\Delta(x))\leq 1$, it means that the complement of the local graph of $x$, $\overline{\Delta(x)}$,  has smallest eigenvalue at least $-2$. Hence, by Theorem \ref{cla2} and the fact that the only non-complete strongly regular graph with smallest eigenvalue bigger than $-2$ is the pentagon, we obtain that the complement of the local graph of $x$, $\overline{\Delta(x)}$, is one of the pentagon, the $3\times 3$ grid, the Clebsch graph or the Schl\"{a}fli graph. This shows that $\Gamma$ has one of the four intersection arrays in the proposition and for each intersection array, there is a unique graph.
\epf

\begin{lemma}\label{ll}
Let $\Gamma$ be a $k$-regular graph with $\nu$ vertices and second largest eigenvalue at most one.
Let $\{A, B\}$ be a partition of the vertex set of $\Gamma$ such that $a := \# A >0$ and $b:= \#B$.
Let $Q=\left[ \begin{array}{cc}
k-\alpha  & \alpha \\
\beta & k-\beta
\end{array} \right]$ be the quotient matrix of the partition $\{A, B\}$. Then $\alpha\geq\frac{(k -1)b}{\nu}$ and, if equality holds, then the partition $\{A, B\}$ is equitable.
\end{lemma}
\pf
By interlacing (Theorem \ref{god}), we obtain that $\alpha+\beta\geq k-1$ as $\theta_1(Q) \leq \theta_1(\Gamma)\leq 1$. Since the size of $A$ is $a$ and the size of $B$ is $b$, we have $\alpha a=\beta b$. These two formulae imply that $\beta\geq \frac{(k-1)a}{\nu}$ and $\alpha\geq\frac{(k -1)b}{\nu}$. Moreover, if equality holds, then the interlacing is tight, and hence $\{A, B\}$ is an equitable partition of $\Gamma$.
\epf

An immediate consequence of Lemma \ref{ll} is the following:

\begin{lemma}\label{lem}
Let $\Gamma$ be a strongly regular graph with prarmeters $(\nu, k, \lambda, \mu)$ such that $\theta_1(\Delta(x))\leq 1$ for any vertex $x$ and let $u$ and $v$ be two fixed vertices at distance 2. Let $Q=\left[ \begin{array}{cc}
\lambda-\alpha  & \alpha \\
\beta & \lambda-\beta
\end{array} \right]$ be the quotient matrix of the partition $\{A, B\}$ of $V(\Delta(v))$, where $A=\Gamma_1(u)\cap \Gamma_1(v)$ and $B=\Gamma_2(u)\cap \Gamma_1(v)$. Then $\beta\geq \frac{(\lambda-1)\mu}{k}$ and $\alpha\geq\frac{(\lambda -1)(k-\mu)}{k}$. Moreover, if equality holds in either of them, then $\{ A, B\}$ is an equitable partition.
\end{lemma}

As another consequence of Lemma \ref{ll}, for regular subgraphs of the complement of the Schl\"{a}fli graph, we obtain the following lemma.

\begin{lemma}\label{schl}
Let $\Gamma$ be the complement of the Schl\"{a}fli graph. Let $\Sigma$ be a $(t+2)$-regular subgraph of $\Gamma$ with $3(t+1)$ vertices, for some $t =0, 1, \ldots, 7$. Then the induced subgraph $\Sigma'$ on $V(\Gamma) \setminus V(\Sigma)$ is a $(9-t)$-regular subgraph (with $24-3t$ vertices).
\end{lemma}

\pf It follows immediately from Lemma \ref{ll}, as $\alpha = t+24$ and $\beta = 9-t$, where $A = V(\Sigma)$ and $B = V(\Gamma) \setminus V(\Sigma)$.
\epf

The following result was shown by B\"{o}inck \cite{boinck}, but for the convenience of the reader we include its proof.

\begin{prop}\label{unique}
There is a unique distance-regular graph $\Gamma$ that is locally the folded 5-cube and with intersection array
$\{ 16, 10, 1; 1, 5, 16\}$.
\end{prop}
\pf
We already have seen the existence of such a graph (Remark \ref{example}).

Now we will show the uniqueness of $\Gamma$.

Fix $x$ a vertex. We will label the vertices of $\Delta(x)$ by the subsets of $\{1, 2,3,4,5\}$ of size at most two, where two subsets $A$ and $B$ are adjacent if $|A \triangle B| = 1$ if one of $A$, $B$ has size at most one, or if $A$ and $B$ have both size two, then $A$ and $B$ are adjacent  if $A \cap B = \emptyset$.  Instead of $\emptyset$ we write $0$, instead of $\{i\}$ we write $i$ and instead of $\{i,j\}$ we write $ij$.

For $y \in \Gamma_2(x)$, let $C(y) := \Gamma(x) \cap \Gamma(y)$ and we define
${\cal C}:= \{ C(y) \mid y \in \Gamma_2(x)\}$.
Then the induced subgraph on $C(y)$ is a pentagon of $\Delta(x)$.

We first give some properties of the set $\cal C$, which are easily checked.
\\
\\
(i) $|C(y) \cap C(z)| \leq 2$ for all $y, z \in \Gamma_2(x)$;
\\
(ii) For each edge $uv$ of $\Delta(x)$, there are exactly four $y \in \Gamma_2(x)$ such that $u, v \in C(y)$;\\
(iii)  For fixed $1 \leq i < j \leq 5$, there are exactly ten $ y \in \Gamma_2(x)$ such that $C(y)$ contains $i, j$;
\\
(iv) For fixed $1 \leq i \leq 5$, there are exactly two $y \in \Gamma_2(x)$ such that $C(y) \cap \{1,2,3,4,5\} = \{i\}$;\\
(v) There are exactly two $y \in \Gamma_2(x)$ such that $C(y) \cap \{1,2,3,4,5\} = \emptyset$;\\
(vi) For $u,v,w \in \Gamma(x)$ such that $u \sim v \sim w$ there is a unique $y \in \Gamma_2(x)$ such $u,v, w \in C(y)$.
\\

We will call pentagons of (iii), pentagons of type 1, pentagons  of (iv), pentagons of type 2 and pentagons of (v), pentagons of type 3.

We will show the set $\cal C$ is unique upto isomorphism.

Let $0', 0'' \in \Gamma_3(0)$. Then $C(0')$ and $C(0'')$ are pentagons of type 3 and $C(0') \cap C(0'') = \emptyset$.

Let $P$ be the subgraph induced on $\Gamma(x) \cap \Gamma_2(0)$. Then $P$ is a Petersen graph.
Moreover, every pentagon of type 1 contains a unique edge of $P$ and every pentagon of type 2 contains a path of length three in $P$.

Now every edge of $P$ can only be the edge of two type 1 pentagons, (by Property (vi)). Also every edge of a pentagon of type 3 is contained in at least one pentagon of type 1. \\
\\
{\bf Claim 1.} There are only two cases, namely \\
(i)  each edge of a type 3 pentagon is contained in exactly two pentagons of type 1, and the other edges are contained in zero pentagons of type 1; and
\\
(ii) each edge of a type 3 pentagon is contained in exactly one pentagon of type 1, and the other edges are contained in two pentagons of type 1.
\\
\\
{\bf Proof of Claim 1:}  For an edge $uv$ of $P$ we define $w(uv)$ as the number of times $uv$ is in a pentagon of type 1. For $uv$ an edge of a type 3 pentagon we have $1 \leq w(uv) \leq 2$ and for the other edges $uv$ we have $0 \leq w(uv) \leq 2$.
We have $\sum_{uv \in E(P)} w(uv) = 20$. If an edge $uv$ of a type 3 pentagon has weight 2, then both edges incident to $uv$, but not lying in a type three pentagon, have weight at most one, which in  turn implies that
$$\sum_{uv \mbox{ not in a pentagon of type 3}} w(uv) \leq 6.$$ This, in turn, implies that there are two incident edge $uv$ and $uw$ of a type both with weight 2. We obtain that the third edge of $P$, containing $u$ has to have weigth 0, and hence
 $$\sum_{uv \mbox{ not in a pentagon of type 3}} w(uv) \leq 4.$$ Continueing in this matter we obtain that all edges whch do not lie in a type 3 pentagon must have weight 0 and the rest weight 2. this shows the claim. \epf
 \\

{\bf Claim 2.} Case (ii) of Claim 1 is not possible.
\\
\\
{\bf Proof of Claim 2:}
 Without loss of generality we may assume that the vertices of $C(0')$ are $12, 34, 15, 23$ and $45$.
 The path $12, 34, 25$ must be in a pentagon of type 2 and hence the fourth vertex of this pentagon contained in $P$ must be $13$ or $14$. In similar fashion the path $35, 12, 34$ must be in a pentagon of type 2 and hence the fourth vertex of this pentagon contained in $P$ must be $24$ or $14$.
 This means that $2 \sim 12 \sim 34 \sim 3$ must lie in a pentagon of type 1. In similar fashion we see that for each edge of a pentagon of type 3, the neighbours of this edge in $\{1,2,3,4,5\}$ are uniquely determined. Now for the edges of $P$ not in  a pentagon of type 3 are edges of two pentagons of type 1. Now, for example, for the edge $12, 35$ either one of these pentagons contains $1, 3 , 12$ and $35$
 and the other one $2, 5, 12$ and $35$. But the edge $35, 24$ lies in a pentagon of type 1 also containing $2$ and $5$ and the edge $35, 14$ lies in a pentagon of type 1 also containing $1$ and $3$.
 But as a pentagon has an odd number of edges, we can not finish the set $\cal C$ in this case.
 \epf
 \\
 \\
{\bf Claim 3:} $\cal C$ is uniquely determined upto isomorphism.
\\
\\
{\bf Proof of Claim 3:}
First we note that the type 2 pentagons are determined by the type 3 pentagons. Now fix a type 2 pentagon $C$. Then there is a unique vertex $u$ of $\Delta(x)$ at distance 2 from this pentagon $C$.
This determines one new pentagon of type 3 with respect to $u$ and at least 4 new pentagons of type 2 with respect to $u$. Continuing in this fashion one easily sees that $\cal C$ is uniquely determined.
\epf
\\
\\
Now the proof is easy to complete. We know the neighbours $y_1, \ldots, y_{10}$ of $0$ in $\Gamma_2(x)$,  and $y_i \sim y_j$ if and only if $i \neq j$ and $C(y_i) \cap C(y_j) \cap \{1,2,3,4,5\} = \emptyset$. As $0$ is any vertex of $\Delta(x)$ we have shown all the edges in $x \cup \Gamma(x) \cup \Gamma_2(x)$. Now an easy induction argument completes the proof.
\epf

\section{Proof of Theorem 1.2}

We will give a proof of Theorem 1.2 in this section. First we will consider the distance-regular graphs with diameter at least three, and later we will consider the strongly regular graphs.

\subsection{Distance-regular graphs with diameter at least three}

Let $\Gamma$ be a distance-regular graph with $\nu$ vertices, diameter $D$ at least three, $a_1 \geq 2$ such that $\theta_1(\Delta(y))\leq 1$ for any vertex $y$ of $\Gamma$.
Clearly $c_2 \geq 2$, otherwise $a_1 \leq 1$. Let $x$ be a fixed vertex of $\Gamma$.
By Proposition \ref{tay} and  Proposition \ref{ter}, we may assume that $\Gamma$ is neither a Taylor graph, nor a Terwilliger graph, and hence we have $c_2 < b_1$, by \cite[Theorem 1.5.5]{bcn}. Since $\overline{\Delta(x)}$ is connected, by Proposition \ref{conn}, we have the following six cases by Theorem \ref{cla1}. \\

Case 1) $\overline{\Delta(x)}$ is the line graph of a $t$-regular graph with $t+\alpha$ vertices where $\alpha \geq 1$ is an integer. \\
Then $k=\frac{1}{2}t(t+\alpha), b_1=2t-2$ and by Lemma \ref{c2}, $c_2\geq \frac{1}{2}t(t+\alpha)-3t+3$.
 Since $b_1 > c_2$, we have $\frac{1}{2}t(t+\alpha)<5t-5$. Hence the only possible pairs of $(t, \alpha)$ are $(4,3)$, ($t =2,4,6$ and $\alpha =2$), and ($2 \leq t \leq 7$ and $\alpha =1$). As $k_2 = kb_1/c_2$ (with $c_2 < b_1$) is an integer and $a_1 \geq 2$, the only possibilities are: $4 \leq t  \leq 6$ and $\alpha =1$, and $t=4$ and $\alpha =2$. By Theorems \ref{Ter} and \ref{Tercla}, we see that the diameter $D$ is three (as for $(\alpha, t) = (1,4)$ we have $k =10$, $b_1 =6$, and hence $D \leq 2\times10/{3+2} = 4$, but equality can not occur).
It is easy to check that the number of vertices of $\Gamma$ is at most 162. Hence, by the tables of \cite[Chapter 14]{bcn}, we see that $\Gamma$ cannot be primitive, and we already assumed $\Gamma$ is not a Taylor graph, so $\Gamma$ must be an antipodal $r$-cover with diameter three, and $r \geq 3$. But this means that $2c_2 \leq b_1$, and the only possible intersection arrays are:
$\{12, 6,1;1, 3, 12\}$, $\{15, 8, 1;1, 4, 15\}$, $\{10,6,1;1,2,10\}$ and $\{10, 6, 1; 1, 3, 10\}$. But there is no distance-regular graph with any of these intersection arrays, see for example \cite[Table 1]{dhks}. \\

Case 2) $\overline{\Delta (x)}$ is the line graph of a bipartite semiregular graph of valency $s$, $t$ with $2 \leq s<t$ (as $a_1 \geq 2$) and, with $\sigma t=\tau s$ vertices ($\sigma \geq s, \tau \geq t)$. \\
Then $k=\sigma t,  b_1=s+t-2$ and by Lemma \ref{c2}, we have  $c_2\geq \sigma t -2s -t +3$. Since $b_1 > c_2$, it follows $t(\sigma -2)<3s-5$.
It follows that $\sigma = s=2.$ But, then $k = 2t$, $b_1 = t$ and $c_2 = t-1$, and hence $t =3$. By Theorems \ref{Ter} and \ref{Tercla}, this is  impossible.\\

Case 3) $\overline{\Delta(x)}$ is a subgraph of $E_7(1)$. Then the possible pairs for $(k, b_1)$ are
$(2t+4, t)$ with $12 \geq t \geq 3$.
By Theorem \ref{Ter}, we obtain $t\geq 7$ and $D=3$. As $\Gamma$ is not Terwilliger we obtain $c_2\geq 2(a_1+1)-k +1+1 =6$. This means that $\Gamma$ has at most $1+28 + 2.28  + 4.28 =197$ vertices.
Hence, by the tables of \cite[Chapter 14]{bcn}, we see that $\Gamma$ cannot be primitive, and we already assumed $\Gamma$ is not a Taylor graph, so $\Gamma$ must be an antipodal $r$-cover with diameter three, and $r \geq 3$. But this means that $2c_2 \leq b_1$, and the only possible intersection array is:
$\{28, 12,1; 1, 6, 28\}$. But there is no distance-regular graph with this intersection array, as by \cite[p.431]{bcn}, its eigenvalues must be integer, but this is not the case.\\

Case 4) $\overline{\Delta(x)}$ is a subgraph of the Schl\"{a}fli graph. Then the possible pairs for $(k, b_1)$ are
$(3(t+1), 2t)$ with $8 \geq t \geq 2$. By Theorem \ref{Ter}, we obtain $D \leq 4$ if $t \geq 5$ and $D=3$ otherwise. Again, with Theorem \ref{Ter} we also see $b_2 - c_2 \leq t-5 \leq 3$. This means that $ \frac{b_2}{c_2} \leq \frac{5}{2}$.
This means that we have one of the following:\\
(a) $D=3$ and $n \leq 784$;\\
(b) $D = 4$ and $n \leq 2134$.
So this means, again by \cite[Chapter 14]{bcn} that $\Gamma$ must be an antipodal $r$-cover with $r \geq 3$. As one of $a_1 = c_2$ and the eigenvalues are integral holds, we obtain that $\Gamma$ has one of the following intersection arrays:
$\{ 27, 16, 1; 1, 4, 27\}$, $\{24, 14, 1; 1, 7, 24\}$, $\{21, 12, 1; 1, 4, 21\}$ and $\{15, 8, 1; 1, 4, 15\}$,
but
no distance-regular graphs exist with intersection array $\{21, 12, 1; 1, 4, 21\}$ and $\{15, 8, 1; 1, 4, 15\}$, as the first one has non-integral multiplicities and the second one does not exist, by \cite[Table 1]{dhks}.  So the two remaining intersection arrays are
$\{ 27, 16, 1; 1, 4, 27\}$, and $\{24, 14, 1; 1, 7, 24\}$.
\\

%Then by  Theorem \ref{Ter} we obtain that possible pairs of $(k, b_1, c_2)$ are $(27,16,\geq 8)$, $(24,14,\geq7)$, $(21,12,\geq6)$, $(18,10\geq5)$, $(15,8,\geq3)$, $(12,6,\geq2)$ and $(9,4,\geq2)$. These imply that $\mid V\Gamma \mid\leq 406$. But by the \cite[Chapter 14]{bcn}, there doesn't exist intersection array. Specially, for the antipodal $r$-cover with $\lambda\neq\mu$, unique possible pairs of $(k, b_1, c_2)$ are $(24,14,7)$ and $(15,8,4)$ as $k=mn$, $\mu=\frac{(m-1)(n+1)}{r}$ and $\lambda=\mu+n-m$ for some integers $-m$ and $n$ which are eigenvalues of $\Gamma$. But since multiplicity of $n$ is less then $k$, these cases can not occur by the Theorem \ref{4.4.4}. \\

Case 5) $\overline{\Delta(x)}$ is a subgraph of Clebsch graph. Then $(k, b_1) = (4t+4 , 3t +1)$, with $1 \leq t \leq 3$. Then by Theorem \ref{Ter} we obtain $D \leq 4$ and $n \leq 1 + 16 + 5.16 + 25. 16 = 497$ if $D =3$ and $n \leq 497 + 125.16 = 2497$ if $D=4$. So again we only need to look at the antipodal $r$-covers with $r \geq 3$, and in similar fashion as in previous case we obtain that $\Gamma$ has $\{ 16, 10, 1; 1, 5, 16\}$ as its intersection array. But then $\Gamma$ is locally the folded 5-cube and there is a unique such distance-regular graph, by Proposition \ref{unique}.
\\

%that possible pairs of $(k, b_1, c_2)$ are $(16,10,\geq 4)$, $(12,7,\geq 3)$ and $(8,4,\geq 2)$. These imply that $\mid V\Gamma \mid\leq 407$. By the \cite[Chapter 14]{bcn}, we obtain possible intersection array $\{16 ,10 ,1 ;1 ,5 ,16 \}$. Specially, for the antipodal $r$-cover with $\lambda\neq\mu$, there doesn't possible intersection array as $k=mn$, $\mu=\frac{(m-1)(n+1)}{r}$ and $\lambda=\mu+n-m$ for some integers $-m$ and $n$ which are eigenvalues of $\Gamma$.\\

Case 6) If $\overline{\triangle (x)}$ is a $K_{m\times 2}$, then $a_1=1$.
\\

So the theorem is shown if the diameter is at least three except that $\Gamma$ can still have two remaining intersection arrays $\{ 27, 16, 1; 1, 4, 27\}$, and $\{24, 14, 1; 1, 7, 24\}$ as its intersection array. As both occur only in Case 4, we see that in the first case $\Gamma$ is locally the complement of the Schl\"{a}fli graph and in the second case locally the subgraph of the complement of Schl\"{a}fli graph in which a triangle is removed (Lemma \ref{schl}). But that means that in both cases the subgraph on the common neighbours of two vertices at distance two has minimal degree at least four, but that means that this subgraph has triangles, which is a contradiction with the fact that the complement of the Schl\"{a}fli graph has no induced $K_{2,1,1}$. This completes the proof of the theorem in case the diameter is at least three. \epf

\subsection{Distance-regular graphs with diameter two }
Let $\Gamma$ be a distance-regular graph with diameter two,  $a_1\geq 2$ and with non-trivial eigenvalues $\theta_1 > \theta_2$ such that $\theta_1(\Delta(x))\leq 1$ for any vertex $x$. We may assume that $\Gamma$ contains a quadrangle and that it is not complete multipartite  and let $x$ be a fixed vertex. By Theorem \ref{cla2} and Proposition \ref{cla3}, we may assume $\theta_1 >1$ and $\theta_2 < -2$. Either $\Gamma$ has only integral eigenvalues or $\Gamma$ has intersection array $\{2t, t; 1, t\}$ with $t\geq 1$ an integer (Lemma \ref{srg} (iv)). For $t \leq 2$, the smallest eigenvalue  $\theta_2$ is at least $-2$. For $t=3$ and $t=4$ there exists a unique graph (see \cite{srg}) namely the Paley graph on, respectively, 13 and 17 vertices. In each case it is easy to check that $\theta_1(\overline{\Delta(x) }) \leq 1$ for all vertices $x$. Also (see \cite{srg}) there does not exist such a graph with intersection array $\{ 21, 10, 1, 10\}$. So we may assume $\theta_1 \geq 2$ and $\theta_2 \leq -3$. \\

For a vertex $x$,  the number $m_x$ will denote the multiplicity of $1$ as an eigenvalue of $\Delta(x)$. If $m_x \geq k_2 +2$ for some vertex $x$, then, by interlacing (Theorem \ref{0}), we see that $\theta_1 = 1$.  So, from now on, we may assume that $m_x$ is at most $1+k_2$ for any vertex $x$. \\

We first will show:\\

{\bf Claim:} $\Gamma$ has possibly one of the following intersection arrays:  $\{ 45, 16,; 1, 24\},$ $\{28,12; 1,16\}$,
$\{27, 16; 1, 6\}$, $\{27, 16; 1, 12\}$, $\{24, 12; 1, 6\}$, $\{21, 12; 1, 6\}$, $\{21, 12; 1, 9\}$, $\{18, 10; 1, 6\}$, $\{15, 8; 1, 6\}$, and $\{12, 6; 1, 6\}$. \\
\\
{\bf Proof of Claim:}

As $\Gamma$ has smallest eigenvalue at most $-3$ and second largest eigenvalue at least two, we find by Theorem \ref{srg}, that $b_1 \geq 6$.

Fix $x$ a vertex of $\Gamma$.
Since $\overline{\Delta(x)}$  is connected for any vertex $x$ (Proposition \ref{conn}), we have the following six cases to consider, by Theorem \ref{cla1}. \\

Case 1) $\overline{\Delta(x)}$ is the line graph of a $t$-regular graph $\Sigma$ with $t+\alpha$ vertices. By looking at the vertex-edge incidence matrix of $\Sigma$ one sees: $\frac{1}{2}t(t+\alpha)-(t+\alpha)\leq m_x$. We also obtain $b_1=2t-2 \geq 6$, $k=\frac{1}{2}t(t+\alpha)$ and, by Lemma \ref{c2},  $c_2\geq k-3t+3$.
Suppose that $k$ is at most $5t-6$. Then $4\leq t\leq 7$.

Suppose that $k$ is at least $5t-5$. Then $m_x \leq1+k_2\leq 1+\frac{2t-2}{k-3t+3}k$, as $c_2\geq k-3t+3$ and $b_1=2t-2 \geq 6 $. This implies $\frac{1}{2}t(t+\alpha)-(t+\alpha)\leq 1+k_2 \leq 1+\frac{2t-2}{k-3t+3}k$. Then $4\leq t\leq10$.

So, in conclusion, we have $4 \leq t \leq 10$.

If $\Gamma$ is a conference graph, then $\frac{t(t+ \alpha)}{2} = k = 2b_1 = 4t-4$, and hence $\Gamma$ has intersection array $\{ 12, 6; 1, 6\}$ (as $t \geq 4$). Now we may assume that $\Gamma$ has integral eigenvalues.
If $t = 4, 5, 6, 8$, then the non-trivial eigenvalues of $\Gamma$ are
$-3$ and $t-2$ (Theorem \ref{srg} (iii)) and hence $c_2 = k-3(t-2)$, by Theorem \ref{srg} (i).
This implies (using $k- t - \alpha \leq k_2 +1$,  and $c_2 = k-3(t-2)$ for $t=4,5,6,8$ and $c_2 \geq k-3t +3$ otherwise,) that for the pair $(t, \alpha)$, we have only the following possibilities: ($t =4$ and $\alpha \leq 5$), ($t=5$ and $\alpha =1, 3$), ($t=6$ and $\alpha =1, 2, 3$), $(t=7$ and $\alpha=1, 3, 5$) and $(t=8, 9, 10$ and $\alpha =1). $
By checking the tables of \cite{srg}, we see that $\Gamma$ has one of the following intersection arrays: $\{ 45, 16; 1, 24\}$, $\{28, 12; 1, 16\}$, $\{15, 8; 1, 6\}$ and $\{12, 6; 1, 6\}$. \\

Case 2 ) $\overline{\Delta(x)}$ is the line graph of a bipartite semiregular graph $\Sigma$ with valencies $s$, $t$ $(2\leq s< t)$ and $\sigma t=\tau s$ edges. Then we obtain $m_x\geq \sigma t-\sigma -\tau$, $b_1=s+t-2$, $k=\sigma t=\tau s$ and, by Lemma \ref{c2}, we have $c_2\geq k-2s-t+3$.
As $b_1 \geq 6$, we obtain $t \geq 5$.

If $\sigma \leq 3$, then $3 \geq \sigma \geq s$, and
hence $s =2,3$. Now $\Gamma$ has a coclique of size at least $t$. As $\theta_1 \geq 2$, we obtain, by
Theorem \ref{srg} (v), that $\Gamma$ has more then $5t -4$ vertices, if $\sigma=2$, and $\Gamma$ has at least $\frac{6t^2}{t+3} + t$ vertices, if $\sigma = 3$. As $c_2 \geq t-1$, if $\sigma =2$, and  $c_2 \geq 2t -3$, if $\sigma =3$, we see that one of the following holds: $\sigma =2$ and $t\leq 7$;
$\sigma =3$ and $t \leq 6$. As $b_1 \geq 6$  we see that the only possibilities for $(s,t,\sigma)$ are $(2, 6, 2)$, $(2, 6, 3)$, $(2, 7, 2)$, $(3, 5, 3)$ and $(3, 6, 3)$.
%Suppose that $k$ is at most $3s+2t-6$. Then $t(\sigma-2)\leq 3s-6$. Since $2\leq s<t$ and $\sigma\geq s$, we obtain $\sigma=2$ and this implies $s=2$. Since $\Gamma$ has coclique of size $t$, $k=2t$, $a_2\leq t$ and $\nu\leq 4t+1$ unless $t=3$, we see that $\theta_1\leq \frac{5}{3}$ by Theorem \ref{srg} (v), contradiction with $\theta_1\geq 2$.

%So, we assume that $k$ is at least $3s+2t-5$.
%If $\Gamma$ is a conference graph then $k = 2b_1 = 2(s+t -2) \leq 4t -6$, so, if $\sigma \geq 4$, then $\Gamma$ has to have integral eigenvalues.

We now consider the case $\sigma \geq 4$. Then $\Gamma$ cannot be a conference graph as $k \neq 2b_1$, so $\Gamma$ has integral eigenvalues.

Suppose that $k$ is at most $3s+2t-6$. Then $t(\sigma-2)\leq 3s-6$. Since $2\leq s<t$ and $\sigma\geq s$, we obtain $\sigma=2$.

So, we assume that $k$ is at least $3s+2t-5$. Then $1+k_2\leq1+\frac{s+t-2}{k-2s-t+3}k\leq 3s+2t-4$ as $\frac{s+t-2}{k-2s-t+3}k$ is decreasing in $k$. As $m_x\leq k_2+1$, we find $\tau s-\sigma-\tau\leq 3s+2t-4$, and hence $(s-2)\tau+5\leq3s+2t$, as $\sigma<\tau$. Now, as $\tau\geq t$, we obtain: $(s-4)t+5\leq3s$.

So, if $s\geq 5$, then for the pair $(s,t)$ we obtain the following possibilities: $s=5$ and $6\leq t\leq 10$. Moreover, if $s=5$, then $3\tau+5\leq 2t+15$, so $6\leq t\leq \tau \leq 10$. This in turn implies $\sigma=5$ and $\tau=t$ as $\sigma t = \tau s$.

For $s=4$ we obtain $2\tau+5\leq 2t+12$ and this implies $\tau\leq t+3$. As $\tau s-\sigma-\tau\leq 3s+2t-4$ and $\sigma t = k = 4 \tau \leq 4t+12$, we obtain $\sigma = s = 4$ and $\tau = t$.

Before we treat $s=2, 3$, we first look at the case $4 \leq s=\sigma\leq 5$, and hence $t=\tau$. Using $1+k_2\leq 1+\frac{s+t-2}{k-2s-t+3}k$ and $k=st$, we obtain $st-s-t\leq 1+\frac{s+t-2}{st-2s-t+3}st$. It is easy to see that there is no solution for $s=5$, for $s=4$, one has $t=5,6$, but $t=5$ is impossible, as then $b_1$ would be a prime.

Now let us return to the case $s=3$ (and $\sigma \geq 4)$). Using $ \tau + 1 \leq 2 \tau - \sigma = k - \tau - \sigma
\leq 1 + k_2 \leq 1 + 3 \tau \frac{t+1}{(\sigma -1)t -3}$, we obtain $4= \sigma$ or $t \leq 6$.
In case $\sigma =4$, in similar fashion as above, one again obtains $t \leq 6$.
As $b_1 = t +1$ is a composite at least 6, we obtain $t = 5$. Now the non-trivial eigenvalues are $-3$ and $2$ and hence $k_2 = 5\sigma  \frac{6}{5 \sigma -6}$ is an integer, but there are no solutions for $\sigma$.

So we are left with the case $s=2$.
Then $b_1 = t$, $k = \sigma t = 2 \tau$.
Then,  as above, $k - \sigma -\tau \leq 1+\frac{s+t-2}{k-2s-t+3}k$,  and hence  we obtain $\sigma =4$ and $t = 6,8$ (as $b_1 = t \geq 6$ and composite).

Summarizing for this case: for the pair $(k, b_1)$ we have the following possibilities:
$(32, 8)$, $(24, 8)$, $(24, 6)$, $(18, 7)$,$(18, 6)$, $(15, 6)$, $(14, 7)$ and $(12, 6)$.

By checking the tables of \cite{srg}, in this case $\Gamma$ has possibly one of the following intersection arrays:
\\
$\{ 12, 6; 1, 6\}$, $\{14, 7; 1, 7\}$,  and $\{15, 6; 1, 9\}$. But  the two intersection arrays $\{ 2t, t; 1, t\} \ \ (t=6,7)$  are impossible in this case, as $\Gamma$ would have a coclique of size $t$ and hence second largest
eigenvalue smaller than two, by Theorem \ref{srg} (v), a contradiction.
And the intersection array $\{15, 6; 1, 9\}$ is impossible, as then the coclique size would be 5, and again this contradicts Theorem \ref{srg} (v).
\\

% As $\sigma t = \tau s$ we also have
% ($s=3= \sigma$ and $4 \leq t \leq 9$), ($s=4= \sigma$ and $5 \leq t \leq 7$), ($s=5= \sigma$ and $t=6$), where we used $\sigma t - \sigma - \tau >\sigma(t-2)$, $\sigma t-\sigma -\tau\leq 3s+2t-5$ and $\sigma t = \tau s$.  If $s=2$, then $

 %So, $s\leq6$. If $s=6$, then $5t-6 \leq \sigma t-\sigma - \frac{\sigma t}{s} \leq 1+\frac{s+t-2}{k-2s-t+4}k\leq 1+ \frac{t+4}{5t-8}6t$ as $s\leq \sigma$. And this implies $t\leq 4$, contradiction. Similarly, if $s=5$, then $t\leq 5$, contradiction. Similarly if $s=4$, then we obtain $t\leq 6$. But there doesn't exist possible intersection array by Theorem \ref{srg}. If $s=3$, then we obtain $t\leq 14$. By Theorem \ref{srg}, we obtain possible intersection arrays $\{15,6;1,9\}$. If $s=2$, then $\sigma(\frac{1}{2}t-1)\leq1+\frac{t}{k-t}k$. For $\sigma\geq 5$, we obtain $5(\frac{1}{2}t-1)\leq 1+\frac{t}{4t}5t$ as $k\geq 5t$. And this implies $t\leq4$, contradiction. Similarly for $\sigma=4$, we obtain $t\leq 7$. But there doesn't exist possible intersection array by Theorem \ref{srg}. For $\sigma=3$, $\Gamma$ has a coclique $C$ of size $t$.
%Since $k=3t$, $a_2\leq t$ and $\mid \Gamma V \mid\leq\frac{9}{2}t+1$, $\theta_1\leq \frac{23}{18}$ by the Claim \ref{claim}. But there doesn't exist possible intersection array by Theorem \ref{srg}. For $\sigma=2$, $\Gamma$ has a coclique $C$ of size $t$. Since $k=2t$, $a_2\leq t$ and $\mid \Gamma V \mid\leq4t+1$, $\theta_1\leq\frac{5}{3}$ as the Claim \ref{claim}. By Theorem \ref{srg}, we obtain possible intersection arrays $\{6,3;1,3\}$, $\{8,4;1,4\}$. \\

Case 3) $\overline{\Delta(x)}$ is a subgraph of $E_7(1)$. Then the possible pairs of $(k, b_1)$ are $(2t+4, t) \ \ (1 \leq t \leq 12)$. As $\Gamma$ is not a conference graph, we find that the eigenvalues of $\Gamma$ must be integral and $\theta_1 \geq 2 $ and $\theta_2 \leq -3$. Hence $b_1 \geq 6$ and a composite. Using $c_2 - k = \theta_1 \theta_2$ and $k_2$ is an integer, we have the following remaining cases:
$\{ 24, 10; 1, 12\}$ and $\{28, 12; 1, 16\}$, but the first case cannot happen, see the tables of \cite{srg}. \\

%By Theorem \ref{srg}, we obtain possible intersection array $\{28, 12; 1, 16\}$. In this case, $\alpha$ in the Lemma \ref{lem} is 5 as $\mu$ of local graph is 10. But, $\alpha\geq\frac{14\cdot 12}{28}=6$, contradiction. \\

Case 4) $\overline{\Delta(x)}$ is a subgraph of the Schl\"{a}fli graph. Then the  possible pairs of $(k, b_1, c_2)$ are $(3t+3,2t)$ \ \ $(1 \leq t \leq 8)$. Then $\Gamma$ has integral eigenvalues. Now in similar fashion  as the previous case, we obtain the following possible intersection arrays: $\{27, 16; 1, 6\}$, $\{27, 16; 1, 12\}$, $\{24, 12; 1, 6\}$, $\{21, 12; 1, 6\}$, $\{21, 12; 1, 9\}$, $\{18, 10; 1, 6\}$, $\{15, 8; 1, 6\}$ and $\{12, 6; 1, 6\}$. \\

Case 5) $\overline{\Delta(x)}$ is a subgraph of the Clebsch graph. Then the possible pairs of $(k, b_1)$ are $(16, 10)$, $(12, 7)$, $(8, 4)$ and $(4, 1)$. Now in similar fashion as in cases 3) and 4), there are no possible intersection arrays.  \\

Case 6) $\overline{\Delta(x)}$ is the $K_{m\times 2}$. Then $a_1=1$. \\

This finishes the proof of the claim.
\epf

To finish the proof of the theorem when the diameter is two, we still need to rule out the intersection arrays: $\{45, 16; 1, 24\}$, $\{28,12; 1, 16\}$,
$\{27, 16; 1, 6\}$, and  $\{24, 12; 1, 6\}$.

When $k=28$, then we see that we are in Case 1) or Case 3). In either case, for any vertex $x$, the local graph  $\Delta(x)$ is the complement of a Chang graph or the complement of $J(8,2)$. This means that $\Delta(x)$ is strongly regular and has $\mu = 10$ and hence for the $\alpha$ in Lemma \ref{lem}, we have $6< \alpha = 6$, a contradiction. In similar fashion, the case $k=45$ is ruled out.

The intersection arrays $\{27, 16; 1, 6\}$, and  $\{24, 12; 1, 6\}$ are ruled out in the same manner as for the intersection arrays $\{ 27, 16, 1; 1, 4, 27\}$, and $\{24, 14, 1; 1, 7, 24\}$ in the diameter three case.
\epf

\section{Proof of Theorem 1.3}

In this section we give a proof of Theorem 1.3. First we recall a result of Hoffman.
Hoffman \cite{hoff} showed:
\begin{theorem}(Cf. \cite[Theorem 3.12.5]{bcn})\label{hoff}
Let $\sigma_k$ be the supremum of the smallest eigenvalues of graphs with minimal valency $k$ and smallest eigenvalue $< -2$. Then $(\sigma_k)_k$ forms a monotone decreasing sequence with limit $-1 - \sqrt{2}$.
\end{theorem}

{\bf Proof of Theorem 1.3:}
Let $\alpha < 1 + \sqrt{2}$.
Let $\Gamma$ be a distance-regular graph with diameter at least three and smallest eigenvalue at most
$-1 - \frac{b_1}{\alpha}$. Let $x$ be a vertex of $\Gamma$. Then $\Delta(x)$ has second largest eigenvalue at most $\alpha-1$. Now the complement of $\Delta(x)$, $\overline{\Delta(x)}$ has smallest eigenvalue at least $-\alpha$ and valency $b_1$. As $D \geq 3$ $b_1 \geq (k+1)/3$ (as $b_1 \geq c_2 \geq 2a_1 +2 -k +1= 2b_1 +k+1$). Theorem \ref{hoff} shows that there exists a $K= K(\alpha)$ such that if $k \geq K$, then $\overline{\Delta(x)}$ has smallest eigenvalue at least $-2$, that is $\Delta(x)$ has second largest eigenvalue at most one. Now by checking the graphs of Theorem \ref{main} we obtain Theorem 1.3.
This completes the proof.
\epf

We end this paper with a remark.

\begin{remark}
(i) A strongly regular graph satisfies $(\theta_1+1)(\theta_2 +1) = - b_1$ and only the conference graphs have non-integral eigenvalues. Therefore Theorem 1.3 is not interesting for $D =2$.
\\
(ii) Question: For given $1 < \beta$ are there only finitely many distance-regular graphs with diameter at least three, $a_1 \geq \beta$ and smallest eigenvalue at most $-1 - \frac{b_1}{\beta}$?
\end{remark}

{\bf Acknowledgments}
\\
We would like to thank Jongyook Park for the careful reading he did. JHK was partially supported by the Basic Science Research Program through the National Research Foundation of Korea(NRF) funded by the Ministry of Education, Science and Technology (grant number 2009-0089826).

\end{document}